\magnification=1200

\def\C{{\bf C}}
\def\D{{\cal D}}
\def\N{{\bf N}}
\def\U{{\cal U}}
\def\hal{{\vrule height 10pt width 4pt depth 0pt}}

\centerline{\bf A counterexample to a conjecture}

\centerline{\bf of Akemann and Anderson}
\medskip

\centerline{Nik Weaver\footnote{*}{Supported by NSF grant DMS-0070634}}
\bigskip
\bigskip
\bigskip

{\narrower{\it
\noindent Akemann and Anderson made a conjecture about ``paving''
projections in finite dimensional matrix algebras which, if true,
would settle the well-known Kadison-Singer problem. We falsify their
conjecture by an explicit seqence of counterexamples.
\bigskip}}
\bigskip

Let $H$ be a separable Hilbert space, let $B(H)$ be the C*-algebra of
bounded linear operators on $H$, and let $\D$ be a maximal abelian
self-adjoint subalgebra (MASA) of $B(H)$. Does every pure state on $\D$
extend uniquely to a pure state on $B(H)$? In [6] Kadison and Singer
answered this question negatively for nonatomic MASAs and speculated
that the answer would also be negative in the atomic case. This has
come to be known as the Kadison-Singer problem.
\medskip

A simple application of the Hahn-Banach theorem shows that every pure state
on $\D$ has at least one extension to a state on $B(H)$, and it is easy
to see that the set of all such extensions is convex and weak* compact.
Hence there exist extreme points in this set, i.e., pure state
extensions. Thus, pure state extensions always exist, and the issue
is uniqueness.
\medskip

Concretely, we can take $H$ to be $l^2(\N)$ and $\D$ to be
$l^\infty(\N)$, regarded as a subalgebra of $B(H)$ via action by
multiplication operators, i.e., diagonal matrices. Since $\D$ is abelian,
the pure states on $\D$ are precisely the homomorphisms from $\D$ into
$\C$, which are well-known to correspond to ultrafilters on $\N$.
(The ultrafilter $\U$ on $\N$ gives rise to the complex homomorphism
$\rho_\U$ defined by $\rho_\U(f) = \lim_\U f$.)
\medskip

It it not too hard to show that if $\U$ is a fixed (trivial)
ultrafilter then $\rho_\U$ has a unique pure state extension.
G.\ Reid [7] showed that there exist free (nontrivial) ultrafilters
$\U$ that also extend uniquely.
\medskip

Let $\pi: B(l^2(\N)) \to l^\infty(\N)$ be the standard conditional
expectation, $\pi(T)(n) = \langle Te_n, e_n \rangle$, where $(e_n)$
is the standard orthonormal basis of $l^2(\N)$. J.\ Anderson [3] showed
that if $\rho$ is any pure state on $l^\infty(\N)$ then $\pi\circ\rho$
is a pure state on $B(l^2(\N))$. Thus, the question of Kadison and
Singer becomes ``Are there ever any other (pure) state extensions of
$\rho$?''
\medskip

The problem can also be rephrased in terms of ``paving'' finite matrices
whose diagonals are zero. The best possible proof of this reduction is
given in [8]. See also [5] for partial positive results and further
references.
\bigskip
\bigskip

\noindent {\bf 1. The Akemann-Anderson conjectures. Idea of the
counterexample}
\bigskip

A different, but related, finite dimensional adaptation of the Kadison-Singer
problem was given by Akemann and Anderson in Chapter 7 of their memoir [1].
They used the following notation. $H = \C^n$ is a finite dimensional complex
Hilbert space and $\D$ is the algebra of diagonal matrices. For any
projection $p \in B(H) = M_n(\C)$, let
$$\delta_p  = {\rm max}\{\langle pe_i, e_i\rangle: 1 \leq i \leq n\}$$
where $\{e_i\}$ is the natural basis of $\C^n$. Recall that a ``symmetry''
is a self-adjoint unitary.
\medskip

Akemann and Anderson state several versions of their conjecture. The primary
version is the following:
\bigskip

\noindent {\bf Conjecture A [1, Conjecture 7.1.1].}
For any projection $p \in M_n(\C)$ there
exists a symmetry $s \in \D$ such that $\|psp\| \leq 2\delta_p$.
\bigskip

This is equivalent to asking whether for any projection $p \in M_n(\C)$
there is a projection $q \in \D$ such that
$${\rm max}(\|qpq\|, \|(1-q)p(1-q)\|) \leq {1\over 2} + \delta_p$$
[1, Proposition 7.4]. This puts the conjecture in the form of a paving
question.
\medskip

If true, Conjecture A
would imply a positive solution to the Kadison-Singer problem
(i.e., unique extension of pure states on atomic MASAs). The grounds for
believing the conjecture include an analogy with known combinatorial
results [2, 4] as well as computer experimentation
which suggests that counterexamples are rare [1, p.\ 80]. Nonetheless,
we will show that this conjecture is false. A weaker version that remains
open, and would still resolve the Kadison-Singer problem, is the following:
\bigskip

\noindent {\bf Conjecture B [1, Conjecture 7.1.3].}
There exist $\gamma, \epsilon > 0$ (independent of $n$) such that for any
projection $p \in M_n(\C)$ with $\delta_p < \gamma$ there exists a symmetry
$s \in \D$ such that $\|psp\| < 1 - \epsilon$.
\bigskip

\noindent I have been unable to falsify this
version of the conjecture. I believe that a counterexample to this version
would probably lead to a negative solution to the Kadison-Singer problem.
\medskip

Since our counterexample is rather computational, it may be helpful to
outline the geometric idea that underlies it. The plan is to produce
a projection $p$ and a unit vector $v$ such that for any symmetry
$s \in \D$ the quantity $\|psp(v)\|$ is larger than $2\delta_p$.
We will take $v$ to be in the range of $p$.
\medskip

Now $s$ is going to be a diagonal matrix whose entries are all $\pm 1$. So,
letting $w_i = p(e_i)$ ($1 \leq i \leq n$), we will have
$$v = p(v) = p\left(\sum \langle v, e_i\rangle e_i\right)
= \sum \langle v, w_i\rangle w_i$$
and
$$psp(v)
= ps(v) = p\left(\sum \pm \langle v, e_i\rangle e_i\right)
= \sum \pm\langle v, w_i\rangle w_i.$$
Thus, the vectors $\langle v, w_i\rangle w_i$ sum to the unit vector $v$,
and the norm of $psp(v)$ for various symmetries $s$
depends on the amount of cancellation which can
be achieved by introducing sign changes in the sum.
\medskip

The idea of the counterexample is to try to arrange $p$ so that
the vectors $w_i' = \langle v, w_i\rangle w_i$ approximate Figure 1. That is,
all but one of them are identical, small, and nearly parallel to $v$, and
the other one has a relatively large component perpendicular to $v$. Then
regardless of how the symmetry $s$ is chosen, the sum which constitutes
$psp(v)$ can be bounded below in norm.
For, if most of the small vectors are given the same sign as the
large vector, then the sum will approximate $\pm v$, which has unit norm;
but if most of the small vectors are given the opposite sign to $v$ then
the sum will have a large component perpendicular to $v$.
\medskip

The actual counterexample is optimized by having not one vector with
a large perpendicular component, but several, each perpendicular in
a different direction. In the notation of the next section, these
terms correspond to the basis vectors $b_i$. The small vectors nearly
parallel to $v$ correspond to the basis vectors $a_i$.
\medskip

In our example the projection $p$ is defined by specifying an orthonormal
basis $\{v_k\}$ of its range, with $v_0$ playing the role of the vector
$v$. The components of the $v_k$ in the $a_i$ and $b_i$ directions are
chosen so that the vectors $\langle v, w_i\rangle w_i$ will have
approximately the form described above. In order to make the $v_k$
orthonormal, and not increase $\delta_p$, we also need to include further
basis vectors $c_{ij}$ and $d_{ij}$.
\bigskip
\bigskip

\noindent {\bf 2. A counterexample to Conjecture A}
\bigskip

Let $m \geq 6$ and let
$H$ be a complex Hilbert space of dimension $2m^3 + 8m^2 + 7m + 2$
(this number is not important). We fix an orthonormal basis for $H$ whose
members fall in four groups, labelled $a$ through $d$. The basis vectors
will be denoted $a_i$ ($1 \leq i \leq m^2$), $b_i$ ($1 \leq i \leq 2m+1$),
$c_{ij}$ ($1 \leq i < j \leq 2m+1$), $d_{ij}$ ($1 \leq i \leq 2m+1$,
$1 \leq j \leq (m+1)^2$). (Thus, the total dimension of $H$ is
$$m^2 + (2m+1) + (2m^2 + m) + (m+1)^2(2m+1) = 2m^3 + 8m^2 + 7m + 2.)$$
$\cal D$ is the MASA of $B(H)$ consisting of those operators which are
diagonal with respect to this basis.
\medskip

We define $p$ to be the orthogonal projection onto the $2m+2$-dimensional
subspace of $H$ which is spanned by the following orthonormal set of
vectors:
$$v_0 = \sum_{i=1}^{m^2} {1\over{m+1}} a_i
+ \sum_{i=1}^{2m+1} {1\over{m+1}} b_i$$
and for $1 \leq i \leq 2m+1$
$$\eqalign{v_i =
&\sum_{j=1}^{m^2} {{-1}\over{m^2(m+1)}} a_j + {1\over{m+1}} b_i\cr
&+ \sum_{j=1}^{i-1} {1\over{m(m+1)}} c_{ji}
+ \sum_{j=i+1}^{2m+1} {{-1}\over{m(m+1)}} c_{ij}
+ \sum_{j=1}^{(m+1)^2} {1\over m} \sqrt{{m-1}\over{m+1}} d_{ij}.\cr}$$

An easy computation shows that these vectors are indeed orthonormal.
\medskip

Next, we show that $\delta_p = 2/(m+1)^2$. We have
$$\eqalign{\|p(a_i)\|^2 &=
\left\|{1\over{m+1}} v_0 + \sum_1^{2m+1} {{-1}\over{m^2(m+1)}} v_i\right\|^2 =
{1\over{(m+1)^2}} + (2m+1)\cdot{1\over{m^4(m+1)^2}}\cr
\|p(b_i)\|^2 &=
\left\|{1\over{m+1}} v_0 + {1\over{m+1}} v_i\right\|^2 = {2\over{(m+1)^2}}\cr
\|p(c_{ij})\|^2 &=
\left\|{{-1}\over{m(m+1)}} v_i + {1\over{m(m+1)}} v_j\right\|^2 =
{2\over{m^2(m+1)^2}}\cr
\|p(d_{ij})\|^2 &=
\left\|{1\over m} \sqrt{{m-1}\over{m+1}} v_i\right\|^2 =
{{m-1}\over{m^2(m+1)}}.\cr}$$
For $m \geq 6$ the largest of these is $2/(m+1)^2$, so this is the value
of $\delta_p$.
\medskip

Now let $s$ be any symmetry in $\D$, so that $s = q-q^\perp$ where $q$
is a projection in $\cal D$. Let
$$\alpha = \#\{i: 1 \leq i \leq m^2, a_i \in {\rm range}(q)\}.$$
By replacing $s$ with $-s$ if necessary, we may suppose $\alpha
\leq m^2/2$. For $1 \leq i \leq m^2$ say $s(a_i) = \epsilon_i a_i$
($\epsilon_i = \pm 1$), and for $1 \leq i \leq 2m+1$ say $s(b_i) =
\epsilon_i' b_i$ ($\epsilon_i' = \pm 1$).
\medskip

We calculate $psp(v_0)$; its norm will give us a lower bound for $\|psp\|$:
$$\eqalign{psp(v_0) =
&\left(\sum_1^{m^2} {1\over{(m+1)^2}} \epsilon_i
+ \sum_1^{2m+1} {1\over{(m+1)^2}} \epsilon_i'\right) v_0\cr
&+\left(\sum_1^{m^2} {{-1}\over{m^2(m+1)^2}} \epsilon_i
+ {1\over{(m+1)^2}}\epsilon_1'\right) v_1\cr
&+ \cdots + \left(\sum_1^{m^2} {{-1}\over{m^2(m+1)^2}} \epsilon_i
+ {1\over{(m+1)^2}} \epsilon_{2m+1}'\right) v_{2m+1}.\cr}$$

Clearly the norm of $psp(v_0)$ exceeds the norms of its projections
onto span$(v_0)$ and span$(v_1, \ldots, v_{2m+1})$. So if $\alpha \leq
m^2/4$ then (by projection onto span$(v_0)$)
$$\eqalign{\|psp(v_0)\|
&\geq \left|\sum_1^{m^2} {1\over{(m+1)^2}} \epsilon_i
+ \sum_1^{2m+1} {1\over{(m+1)^2}} \epsilon_i'\right|\cr
&= \left|{{2\alpha-m^2}\over{(m+1)^2}}
+ \sum_1^{2m+1} {1\over{(m+1)^2}} \epsilon_i'\right|\cr
&\geq {{m^2}\over 2}\cdot{1\over{(m+1)^2}} - (2m+1)\cdot{1\over{(m+1)^2}}\cr
&\geq {{m^2 - 4m - 2}\over 4} \delta_p;\cr}$$
while if $\alpha \geq m^2/4$ then (by projection onto span$(v_1, \ldots,
v_{2m+1})$)
$$\eqalign{\|psp(v_0)\|
&\geq \sqrt{2m+1} \cdot \left|\sum_1^{m^2} {{-1}\over{m^2(m+1)^2}} \epsilon_i
\pm {1\over{(m+1)^2}}\right|\cr
&= \sqrt{2m+1} \cdot \left|{{m^2-2\alpha}\over{m^2(m+1)^2}}
\pm {1\over{(m+1)^2}}\right|\cr
&\geq \sqrt{2m+1} \cdot {{2\alpha}\over{m^2(m+1)^2}}\cr
&\geq {{\sqrt{2m+1}}\over 4} \cdot \delta_p.\cr}$$
We conclude that for any symmetry $s$ in $\cal D$,
$$\eqalign{\|psp\|
&\geq {\delta_p\over 4}\cdot{\rm min}(m^2-4m-2, \sqrt{2m+1})\cr
&= \delta_p\cdot\sqrt{m + 1/2}\cr.}$$
\bigskip
\bigskip

\noindent {\bf 3. A limitation of the method}
\bigskip

The counterexample presented in Section 2 is, in a sense, stronger than
needed. To falsify the Akemann-Anderson conjecture, we had to place a
lower bound on
$\|psp\|$ for any symmetry $s \in \D$. That is, for any
$s$, we had to find a unit vector $v$ such that $\|psp(v)\|$ was
sufficiently large. In fact, our counterexample was such that the
same vector $v_0$ could play this role for all symmetries.
\medskip

In pursuing stronger counterexamples, it is worth observing that
there is a limitation to this approach. Namely, for any single
vector $v$ it is always possible to find a symmetry $s$ such that
$\|psp(v)\|$ is on the order of $\delta_p^{1/2}$. This means that
``single vector'' counterexamples cannot falsify Conjecture B.
In detail, we have the following result. We retain the notation of Section 1.
\bigskip

\noindent {\bf Theorem 1.} {\it Let $p \in M_n(\C)$ be a projection. For
any unit vector $v \in \C^n$ there exists a symmetry $s \in \cal D$ such that
$$\|psp(v)\| \leq \sqrt{2\delta_p+3\delta_p^2}.$$}
\bigskip

The proof of the theorem requires the following lemmas.
\bigskip

\noindent {\bf Lemma 2.} {\it Let $(v_i)$ be a sequence of vectors in a
Hilbert space and let $(w_i)$ be the sequence of partial sums,
$w_i = v_1 + \cdots + v_i$. Suppose that
$${\rm Re}\langle v_{i+1},w_i\rangle \leq 0$$
for all $i$. Then
$$\|w_i\|^2 \leq \|v_1\|^2 + \cdots + \|v_i\|^2$$
for all $i$.}
\medskip

\noindent {\it Proof.} By induction:
$$\eqalign{\|w_{i+1}\|^2
&= \|w_i\|^2 + \|v_{i+1}\|^2 + 2{\rm Re}\langle v_{i+1},w_i\rangle\cr
&\leq \|w_i\|^2 + \|v_{i+1}\|^2\cr}$$
and $\|w_1\|^2 = \|v_1\|^2$.\hfill\hal
\bigskip

\noindent {\bf Lemma 3.} {\it Let $(v_i)$ be a finite sequence of vectors
such that $\sum v_i = 0$. Then the sequence can be rearranged so as to
satisfy the hypothesis of Lemma 2.}
\medskip

\noindent {\it Proof.} Take $v_1$ for the first vector in the rearrangement.
Then since $\sum_{i\neq 1} v_i = -v_1$, and hence
$$\left\langle \sum_{i\neq 1} v_i, v_1\right\rangle \leq 0,$$
we must have
$${\rm Re}\langle v_i, v_1\rangle \leq 0$$
for some $i \neq 1$. Take this $v_i$ for the second vector in the
rearrangement, and continue inductively.\hfill\hal
\bigskip

\noindent {\it Proof of Theorem.} It will suffice to consider the
case that $v \in$ range$(p)$ (otherwise replace $v$ by $p(v)$).
Let $\{q_i\}$ be the minimal projections in $\cal D$,
$q_i$ = projection onto span$(e_i)$. Let $p_1$ be the projection onto
span$(v)$ and $p_2 = p - p_1$. For $1 \leq i \leq n$ let $\alpha_i =
\|p_1(e_i)\|$,
$\beta_i = \|p_2(e_i)\|$, $x_i = p_1q_ip(v)$, and $y_i = p_2q_ip(v)$.
Note that $x_i + y_i = pq_ip(v)$; $\|x_i\| = \alpha_i^2$ and
$\|y_i\| = \alpha_i\beta_i$; and
$$\delta_p = {\rm max}\{\|p(e_i)\|^2\} = {\rm max}\{\alpha_i^2 + \beta_i^2\}.$$

Now $\sum pq_ip(v) = v$, so $\sum x_i = v$ and $\sum y_i = 0$. By Lemma
3 we may (rearranging the indices) suppose that the sequence $(y_i)$
satisfies the hypothesis of Lemma 2. Since $\sum x_i = v$ and $x_i
= \alpha_i^2 v$ it follows that $\sum \alpha_i^2 = 1$; also as
$\alpha_i^2 \leq \delta_p$ for all $i$, we can
find an integer $k$ such that $|1/2 - \sum_1^k \alpha_i^2| \leq \delta_p/2$.
\medskip

Let $q = q_1 + \cdots + q_k$ and let
$s = q - q^\perp = 2q - 1$.  Then by the choice of $k$ we have
$$\|p_1sp(v)\|  = \|p_1q(v) - p_1q^\perp(v)\|
= \left|\sum_1^k \alpha_i^2 - \sum_{k+1}^n \alpha_i^2\right|\leq \delta_p,$$
and by Lemma 2,
$$\eqalign{\|p_2sp(v)\|^2 &= \|p_2(2q-1)p(v)\| =
4\|p_2qp(v)\|^2 = 4\left\|\sum_1^k y_i\right\|^2\cr
&\leq 4\sum_1^k \|y_i\|^2
= 4\sum_1^k \alpha_i^2\beta_i^2 \leq 4\delta_p\sum_1^k \alpha_i^2\cr
& \leq 4\delta_p(1/2 + \delta_p/2) = 2\delta_p + 2\delta_p^2.\cr}$$
Thus
$$\|psp(v)\|^2 \leq \delta_p^2 + 2\delta_p + 2\delta_p^2
= 2\delta_p + 3\delta_p^2,$$
as desired.\hfill\hal
\bigskip

In light of this limitation, significant strengthening of the counterexample
will require that the vector $v$ be allowed to vary. The simplest way to
achieve this would be to fix two vectors, $v_0$ and $v_1$, at the outset,
and arrange matters so that no matter how $s$ is chosen either
$\|psp(v_0)\|$ or $\|psp(v_1)\|$ is larger than some desired lower bound
(the dichotomy depending, perhaps, on a parameter similar to $\alpha$ as
in the example of
Section 2). This would be sufficient to bound $\|psp\|$ below, and it
would allow greater flexibility in constructing examples. Using this method
I have been able to construct a series of examples slightly improving
those of Section 2, for which $\|psp\|$ is bounded below by a constant
times $\delta_p^{1/2}$.

\bigskip
\bigskip

[1] C.\ A.\ Akemann and J.\ Anderson, {\it Lyapunov Theorems for
Operator Algebras}, {\it Mem.\ Amer.\ Math.\ Soc.\ \bf 94} no.\ 458 (1991).
\medskip

[2] {---------}, The continuous Beck-Fiala theorem is optimal, {\it
Discrete Math.\ \bf 146} (1995), 1-9.
\medskip

[3] J.\ Anderson, Extreme points in sets of positive linear maps on
$B(H)$, {\it J.\ Funct.\ Anal.\ \bf 31} (1979), 195-217
\medskip

[4] J.\ Beck and T.\ Fiala, ``Integer-making'' theorems, {\it
Discrete Appl.\ Math.\ \bf 3} (1981), 1-8.
\medskip

[5] J.\ Bourgain and L.\ Tzafriri, On a problem of Kadison and Singer,
{\it J.\ Reine Angew.\ Math.\ \bf 420} (1991), 1-43.
\medskip

[6] R.\ V.\ Kadison and I.\ M.\ Singer, Extensions of pure states,
{\it Amer.\ J.\ Math.\ \bf 81} (1959), 383-400.
\medskip

[7] G.\ A.\ Reid, On the Calkin representation, {\it Proc.\ London Math.\
Soc.\ \bf 23} (1971), 547-564.
\medskip

[8] B.\ Tanbay,  Pure state extensions and compressibility of the
$l_1$-algebra, {\it Proc.\ Amer.\ Math.\ Soc.\ \bf 113} (1991), 707-713.
\bigskip
\bigskip

\noindent Math Dept.

\noindent Washington University

\noindent St.\ Louis, MO 63130

\noindent nweaver@math.wustl.edu
\end